\documentclass{amsart}
\usepackage{amssymb}
\usepackage{amsfonts}

\newtheorem{theorem}{Theorem}
\theoremstyle{plain}

\newtheorem{conjecture}{Conjecture}

\newtheorem{lemma}{Lemma}

\newtheorem{proposition}{Proposition}
\newtheorem{remark}{Remark}

\numberwithin{equation}{section}
\input{tcilatex}

\begin{document}
\title[Uniqueness results for positive harmonic functions on $\overline{%
\mathbb{B}^{n}}$]{Uniqueness results for positive harmonic functions on $%
\overline{\mathbb{B}^{n}}$ satisfying a nonlinear boundary condition}
\author{Qianqiao Guo}
\address{Department of Applied Mathematics, Northwestern Polytechnical
University, Xi'an, Shaanxi 710129, China}
\email{gqianqiao@nwpu.edu.cn}
\author{Xiaodong Wang}
\address{Department of Mathematics, Michigan State University, East Lansing,
MI 48824}
\email{xwang@math.msu.edu}
\maketitle

\section{\protect\bigskip Introduction}

From the Hardy-Littlewood-Sobolev inequality on $\mathbb{S}^{n-1}$ with
sharp constant Beckner \cite{B} derived the following family of inequalities
on $\mathbb{B}^{n}$

\begin{theorem}
Let $u\in C^{\infty }\left( \overline{\mathbb{B}^{n}}\right) $. Then 
\begin{equation}
c_{n}{}^{\left( q-1\right) /\left( q+1\right) }\left( \int_{\mathbb{S}%
^{n-1}}\left\vert u(\xi )\right\vert ^{q+1}d\sigma \left( \xi \right)
\right) ^{2/\left( q+1\right) }\leq \left( q-1\right) \int_{\mathbb{B}%
^{n}}\left\vert \nabla u(x)\right\vert ^{2}dx+\int_{\mathbb{S}%
^{n-1}}\left\vert u(\xi )\right\vert ^{2}d\sigma \left( \xi \right) ,
\label{ineq}
\end{equation}
for $1<q<\infty $ if $n=2$ and $1<q\leq n/\left( n-2\right) $ if $n\geq 3$,
where $c_{n}=2\pi ^{n/2}/\Gamma \left( n/2\right) =\left\vert \mathbb{S}%
^{n-1}\right\vert $ and $d\sigma $ is the standard volume form on $\mathbb{S}%
^{n-1}$.
\end{theorem}

\bigskip The critical case $q=n/\left( n-2\right) $ was also proved by
Escobar \cite{E1, E2} by a different method. By Lions \cite{L} the following 
\begin{equation*}
\inf \frac{\int_{\mathbb{B}^{n}}\left\vert \nabla u\right\vert ^{2}+\frac{n-2%
}{2}\int_{\mathbb{S}^{n-1}}u^{2}}{\left( \int_{\mathbb{S}^{n-1}}\left\vert
u\right\vert ^{2\left( n-1\right) /\left( n-2\right) }\right) ^{\left(
n-2\right) /\left( n-1\right) }}
\end{equation*}%
is achieved by a \textbf{positive} function $u$ satisfying the following%
\begin{equation*}
\begin{array}{ccc}
\Delta u=0 & \text{on} & \mathbb{B}^{n}, \\ 
\frac{\partial u}{\partial \nu }+\frac{n-2}{2}u=u^{\left( n-1\right) /\left(
n-2\right) } & \text{on} & \mathbb{S}^{n-1}.%
\end{array}%
\end{equation*}%
Escobar \cite{E2} then classified all positive solutions of the above
equation using an integral method and hence proved the inequality (\ref{ineq}%
) for $q=n/\left( n-2\right) $. The inequality for $1<q<n/\left( n-2\right) $
would also follow in the same way from the following

\begin{conjecture}
If $u\in C^{\infty }\left( \overline{\mathbb{B}^{n}}\right) $ is positive
and satisfies the following equation%
\begin{equation}
\begin{array}{ccc}
\Delta u=0 & \text{on} & \mathbb{B}^{n}, \\ 
\frac{\partial u}{\partial \nu }+au=u^{q} & \text{on} & \mathbb{S}^{n-1},%
\end{array}
\label{E}
\end{equation}
then $u$ is constant, provided $1<q<n/\left( n-2\right) $ and $0<a\leq
1/\left( q-1\right) $.
\end{conjecture}

In fact, when $q<n/\left( n-2\right) $ the trace operator $H^{1}\left( 
\mathbb{B}^{n}\right) \rightarrow L^{q}\left( \mathbb{S}^{n-1}\right) $ is
compact. It follows trivially that 
\begin{equation*}
\inf \frac{\int_{\mathbb{B}^{n}}\left\vert \nabla u\right\vert ^{2}+a\int_{%
\mathbb{S}^{n-1}}u^{2}}{\left( \int_{\mathbb{S}^{n-1}}\left\vert
u\right\vert ^{q+1}\right) ^{2/\left( q+1\right) }}
\end{equation*}%
is achieved by a positive function $u$ satisfying (\ref{E}). The conjecture,
if true, then implies that $u$ is constant if $a\leq 1/\left( q-1\right) $
and hence the inequality (\ref{ineq}). By taking limit $q\nearrow n/\left(
n-2\right) $ one would also obtain the critical case $q=n/\left( n-2\right) $%
. \ In \cite{W2} the conjecture is even formulated in a more general
context, namely the same uniqueness result should be true on any compact
Riemannian manifold with nonnegative Ricci and with principal curvature $%
\geq 1$ on the boundary. In its precise form, the conjecture states

\begin{conjecture}
(\cite{W2}) Let $\left( M^{n},g\right) $ be a smooth compact Riemannian
manifold with $Ric\geq 0$ and $\Pi \geq 1$ on its boundary $\Sigma $. Let $%
u\in C^{\infty }\left( M\right) $ be a positive solution to the following
equation%
\begin{equation}
\begin{array}{ccc}
\Delta u=0 & \text{on} & M, \\ 
\frac{\partial u}{\partial \nu }+au=u^{q} & \text{on} & \Sigma ,%
\end{array}
\label{Eq}
\end{equation}%
where the parameters $\lambda $ and $q$ are always assumed to satisfy $a>0$
and $1<q\leq \frac{n}{n-2}$. If $a\leq \frac{1}{q-1}$, then $u$ must be
constant unless $q=\frac{n}{n-2}$, $M$ is isometric to $\overline{\mathbb{B}%
^{n}}\subset \mathbb{R}^{n}$ and $u$ corresponds to 
\begin{equation*}
u\left( x\right) =\left[ \frac{2}{n-2}\frac{1-\left\vert \xi \right\vert ^{2}%
}{1+\left\vert \xi \right\vert ^{2}\left\vert x\right\vert ^{2}-2x\cdot \xi }%
\right] ^{\left( n-2\right) /2}
\end{equation*}%
for some $\xi \in \mathbb{B}^{n}$.
\end{conjecture}

This conjecture, if true, would yield Beckner-type inequalities on such
manfiolds. We refer to \cite{W2} for further discussion.

Conjecture 2 has been completely confirmed in dimension 2 in \cite{GHW}.

\begin{theorem}
Let $(\Sigma ,g)$ be a compact surface with Gaussian curvature $K\geq 0$ and
geodesic curvature $\kappa \geq 1$ on the boundary. The only positive
solution to the followng equation 
\begin{equation*}
\begin{array}{ccc}
\Delta u=0 & \text{on} & \Sigma , \\ 
\frac{\partial u}{\partial \nu }+au=u^{q} & \text{on} & \partial \Sigma ,%
\end{array}%
\end{equation*}%
where $a>0$ and $q>1$, is constant, provided $a\left( q-1\right) \leq 1$.
\end{theorem}

\bigskip In higher dimensions the following partial result is proved in \cite%
{GHW}.

\begin{theorem}
Let $\left( M^{n},g\right) $ be a smooth compact Riemannian manifold with
nonnegative sectional curvature and the second fundamental form of the
boundary $\Pi \geq 1$. Then the only positive solution to (\ref{Eq}) is
constant if $\lambda \leq \frac{1}{q-1}$, provided $3\leq n\leq 8$ and $%
1<q\leq \frac{4n}{5n-9}$.
\end{theorem}

In this short note we study the problem on the model space $\overline{%
\mathbb{B}^{n}},n\geq 3$. The main result is the following partial result in
dimensions $n\geq 3$.

\begin{theorem}
\label{main}If $u\in C^{\infty }\left( \overline{\mathbb{B}^{n}}\right) $ is
positive and satisfies the following equation (\ref{E}), then $u$ is
constant, provided $1<q<n/\left( n-2\right) $ and $0<a\leq \left( n-2\right)
/2$.
\end{theorem}

The paper is organized as follows. In section 2 we transform the equation (%
\ref{E}) to a new equation on the upper half space $\overline{\mathbb{R}%
_{+}^{n}}=\left\{ x_{n}\geq 0\right\} $. We then study the new equation by
the method of moving planes and prove that the solution is axially symmetric
with respect to the $x_{n}$-axis. In section 3 we go back to the ball and
finish the proof of Theorem \ref{main} by an integral identiy.

\textbf{Aknowlegment. }We want to thank Fengbo Hang and Meijun Zhu for
useful discussions. The 2nd author is partially supported by Simons
Foundation Collaboration Grant for Mathematicians \#312820.

\section{Analysis on $\overline{\mathbb{R}_{+}^{n}}$}

In this section we prove the following partial result.

\begin{proposition}
\bigskip Suppose $u\in C^{\infty }\left( \overline{\mathbb{B}^{n}}\right) $
is positive and satisfies the equation (\ref{E}) with $1<q<n/\left(
n-2\right) $ and $0<a\leq \left( n-2\right) /2$. If $\xi \in \mathbb{S}%
^{n-1} $ is a critical point of \ $u|_{\mathbb{S}^{n-1}}$, then $u$ is axial
symmetric w.r.t. the line through the origin and $\xi $.
\end{proposition}

Without loss of generality, we assume that the north pole $e_{n}$ is a
critical point of $u|_{\mathbb{S}^{n-1}}$. The inverse of the steregraphic
projection $\Psi :\overline{\mathbb{R}_{+}^{n}}\rightarrow \overline{\mathbb{%
B}^{n}}\backslash \left\{ e_{n}\right\} $ is given by%
\begin{equation*}
\Psi \left( x\right) =\left( \frac{2x_{1}}{1+\left\vert x\right\vert
^{2}+2x_{n}},\cdots ,\frac{2x_{n-1}}{1+\left\vert x\right\vert ^{2}+2x_{n}},%
\frac{-1+\left\vert x\right\vert ^{2}}{1+\left\vert x\right\vert ^{2}+2x_{n}}%
\right) .
\end{equation*}%
Let $v\left( x\right) =u\circ \Psi \left( x\right) \left( \frac{2}{%
1+\left\vert x\right\vert ^{2}+2x_{n}}\right) ^{\left( n-2\right) /2}$. Then 
$v$ satisfies the following equation on $\overline{\mathbb{R}_{+}^{n}}$ 
\begin{equation}
\begin{array}{ccc}
\Delta v=0 & \text{on} & \mathbb{R}_{+}^{n}, \\ 
-\frac{\partial v}{\partial x_{n}}=\alpha \left( \frac{2}{1+\left\vert
x\right\vert ^{2}}\right) ^{2}v+\left( \frac{2}{1+\left\vert x\right\vert
^{2}}\right) ^{\beta }v^{q} & \text{on} & \mathbb{R}^{n-1},%
\end{array}
\label{ER}
\end{equation}%
where 
\begin{equation*}
\alpha =\frac{n-2}{2}-a,\beta =\frac{n-q\left( n-2\right) }{2}.
\end{equation*}%
By our assumption $\alpha $ and $\beta $ are both nonnegative. As $e_{n}$ is
a critical point of $u|_{\mathbb{S}^{n-1}}$, using the Taylor expansion of $%
u $ at $e_{n}$ we have, as $x\rightarrow \infty $, 
\begin{eqnarray*}
v\left( x\right) &=&c_{0}\left\vert x\right\vert ^{2-n}\left( 1+c_{1}\frac{%
x_{n}}{\left\vert x\right\vert ^{2}}+O\left( \left\vert x\right\vert
^{-2}\right) \right) \text{ }, \\
\frac{\partial v}{\partial x_{i}} &=&c_{0}\left( -\frac{\left( n-2\right)
x_{i}}{\left\vert x\right\vert ^{n}}-\frac{nc_{1}x_{n}x_{i}}{\left\vert
x\right\vert ^{n+2}}+O\left( \frac{1}{\left\vert x\right\vert ^{n+1}}\right)
\right) ,i=1,\cdots ,n-1,
\end{eqnarray*}%
where $c_{0}$ is a positive constant. We will prove that $v$ is axially
symmetric with respect to the $x_{n}$-axis by the method of moving planes.
We will follow the approach in the classic work of Gidas, Ni and Nirenberg 
\cite{GNN}.

\begin{remark}
When $\alpha =0,\beta =0$ the equation (\ref{ER}) reduces to the following 
\begin{equation*}
\begin{array}{ccc}
\Delta v=0 & \text{on} & \mathbb{R}_{+}^{n}, \\ 
-\frac{\partial v}{\partial x_{n}}=v^{q} & \text{on} & \mathbb{R}^{n-1}.%
\end{array}%
\end{equation*}%
This equation is invariant under translations $x\rightarrow x+\eta $ when $%
\eta _{n}=0$. When $q=n/\left( n-2\right) $ it is further invariant under
Mobius transformations and all postive soluitons were classified by Li and
Zhu \cite{LZ} using the more powerful moving sphere method. Ou \cite{Ou}
studied the case $q<n/\left( n-2\right) $ by the method of moving planes and
proved that there is no positive solution.
\end{remark}

Since the equation (\ref{ER}) is invariant under rotation about the $x_{n}$%
-axis, it suffices to show that $v$ is even in $x_{1}$. For $\lambda \in 
\mathbb{R}$ we define $x^{\lambda }=\left( 2\lambda -x_{1},x_{2,}\cdots
,x_{n}\right) $ and $\Sigma _{\lambda }=\left\{ x\in \mathbb{R}%
^{n}:x_{1}\leq \lambda \right\} $. Let $\Lambda =\left\{ \lambda >0:v\left(
x\right) \geq v\left( x^{\lambda }\right) \text{ on }\Sigma _{\lambda
}\right\} $.

\begin{lemma}
\label{start}For any $\lambda _{0}>0$ there exists $R_{0}>0$ s.t. for all $%
\lambda \geq \lambda _{0},x\in \Sigma _{\lambda }$ with $\left\vert
x\right\vert \geq R_{0}$%
\begin{equation*}
v\left( x\right) \geq v\left( x^{\lambda }\right) .
\end{equation*}
\end{lemma}

By this lemma, it is obvious that $\Lambda $ contains all sufficiently large 
$\lambda $. To prove the lemma, we assume $c_{0}=1$ without loss of
generality in the following proof. First observe that we always have $%
\left\vert x\right\vert \leq \left\vert x^{\lambda }\right\vert $ when $%
\lambda >0$ and $x\in \Sigma _{\lambda }$. We consider several cases.

\textbf{Case 1}. $\left\vert x^{\lambda }\right\vert \geq 2\left\vert
x\right\vert $.

Then 
\begin{eqnarray*}
v\left( x\right) -v\left( x^{\lambda }\right) &=&\left\vert x\right\vert
^{2-n}-\left\vert x^{\lambda }\right\vert ^{2-n}+O\left( \left\vert
x\right\vert ^{1-n}\right) \\
&=&\left\vert x\right\vert ^{2-n}\left( 1-\left( \frac{\left\vert
x\right\vert }{\left\vert x^{\lambda }\right\vert }\right) ^{n-2}+O\left(
\left\vert x\right\vert ^{-1}\right) \right) \\
&\geq &\left\vert x\right\vert ^{2-n}\left( 1-2^{-\left( n-2\right)
}+O\left( \left\vert x\right\vert ^{-1}\right) \right) \\
&\geq &0,\text{ if }\left\vert x\right\vert \text{ is sufficiently large.}
\end{eqnarray*}

\textbf{Case 2}. $\left\vert x^{\lambda }\right\vert \leq 2\left\vert
x\right\vert $, but $\left\vert x^{\prime }\right\vert \leq \left\vert
x\right\vert /\sqrt{2}$. Here we write $x=\left( x_{1},x^{\prime }\right) $.

We have, with $y_{1}=2\lambda -x_{1}$ 
\begin{eqnarray*}
v\left( x\right) -v\left( x^{\lambda }\right) &=&\int_{x_{1}}^{y_{1}}-\frac{%
\partial u}{\partial x_{1}}\left( t,x^{\prime }\right) dt \\
&=&\int_{x_{1}}^{y_{1}}\left[ \frac{\left( n-2\right) t}{\left(
t^{2}+\left\vert x^{\prime }\right\vert ^{2}\right) ^{n/2}}+\frac{%
nc_{1}x_{n}t}{\left( t^{2}+\left\vert x^{\prime }\right\vert ^{2}\right)
^{\left( n+2\right) /2}}\right] dt+O\left( \frac{y_{1}-x_{1}}{\left\vert
x\right\vert ^{n+1}}\right) \\
&=&\frac{1}{\left\vert x\right\vert ^{n-2}}-\frac{1}{\left\vert x^{\lambda
}\right\vert ^{n-2}}+\widetilde{c}x_{n}\left( \frac{1}{\left\vert
x\right\vert ^{n}}-\frac{1}{\left\vert x^{\lambda }\right\vert ^{n}}\right)
+O\left( \frac{y_{1}-x_{1}}{\left\vert x\right\vert ^{n+1}}\right) \\
&=&\frac{\left( n-2\right) \left( \left\vert x^{\lambda }\right\vert
-\left\vert x\right\vert \right) }{r^{n-1}}+\frac{n\widetilde{c}x_{n}\left(
\left\vert x^{\lambda }\right\vert -\left\vert x\right\vert \right) }{s^{n+1}%
}+O\left( \frac{y_{1}-x_{1}}{\left\vert x\right\vert ^{n+1}}\right) ,
\end{eqnarray*}%
where in the last step we used the mean value theorem to get $r,s\in \left(
\left\vert x\right\vert ,2\left\vert x\right\vert \right) $. On the other
hand, as $\left\vert x^{\lambda }\right\vert ^{2}-\left\vert x\right\vert
^{2}=y_{1}^{2}-x_{1}^{2}=2\lambda \left( y_{1}-x_{1}\right) $, we have%
\begin{equation*}
\left\vert y_{1}-x_{1}\right\vert \leq \frac{\left\vert x^{\lambda
}\right\vert ^{2}-\left\vert x\right\vert ^{2}}{2\lambda _{0}}\leq \frac{%
2\left\vert x\right\vert }{\lambda _{0}}\left( \left\vert x^{\lambda
}\right\vert -\left\vert x\right\vert \right) .
\end{equation*}%
Therefore%
\begin{eqnarray*}
v\left( x\right) -v\left( x^{\lambda }\right) &\geq &\frac{\left( n-2\right)
\left( \left\vert x^{\lambda }\right\vert -\left\vert x\right\vert \right) }{%
2^{n-1}\left\vert x\right\vert ^{n-1}}+O\left( \frac{\left\vert x^{\lambda
}\right\vert -\left\vert x\right\vert }{\left\vert x\right\vert ^{n}}\right)
+O\left( \frac{\left\vert x^{\lambda }\right\vert -\left\vert x\right\vert }{%
\lambda _{0}\left\vert x\right\vert ^{n}}\right) \\
&=&\frac{\left\vert x^{\lambda }\right\vert -\left\vert x\right\vert }{%
\left\vert x\right\vert ^{n-1}}\left[ \frac{\left( n-2\right) }{2^{n-1}}%
+O\left( \frac{1+\lambda _{0}^{-1}}{\left\vert x\right\vert }\right) \right]
\\
&\geq &0,\text{ if }\left\vert x\right\vert \text{ is sufficiently large,
depending on }\lambda _{0}\text{.}
\end{eqnarray*}

\textbf{Case 3}. $\left\vert x^{\lambda }\right\vert \leq 2\left\vert
x\right\vert ,\left\vert x^{\prime }\right\vert \leq \left\vert x\right\vert
/\sqrt{2}$.

It follows $\left\vert x_{1}\right\vert \geq \left\vert x\right\vert /\sqrt{2%
}$. If $x_{1}\geq \left\vert x\right\vert /\sqrt{2}$, by the mean value
theorem, there is $s\in \left( x_{1},y_{1}\right) $%
\begin{eqnarray*}
v\left( x\right) -v\left( x^{\lambda }\right) &=&\frac{\partial u}{\partial
x_{1}}\left( s,x^{\prime }\right) \left( x_{1}-y_{1}\right) \\
&=&\left[ \frac{\left( n-2\right) s}{\left( s^{2}+\left\vert x^{\prime
}\right\vert ^{2}\right) ^{n/2}}+O\left( \frac{1}{\left( s^{2}+\left\vert
x^{\prime }\right\vert ^{2}\right) ^{n/2}}\right) \right] \left(
y_{1}-x_{1}\right) \\
&\geq &\left[ \frac{\left( n-2\right) }{\sqrt{2}\left\vert x\right\vert
^{n-1}}+O\left( \frac{1}{\left\vert x\right\vert ^{n}}\right) \right] \left(
y_{1}-x_{1}\right) \\
&\geq &0,\text{ if }\left\vert x\right\vert \text{ is sufficiently large.}
\end{eqnarray*}%
If $x_{1}\leq -\left\vert x\right\vert /\sqrt{2}$, we let $\overline{x}%
=\left( -x_{1},x^{\prime }\right) $. Then By the asymptotic expansion 
\begin{equation*}
v\left( x\right) -v\left( \overline{x}\right) =O\left( \frac{1}{\left\vert
x\right\vert ^{n}}\right) .
\end{equation*}%
By the same type of analysis as in Case 2, we have%
\begin{eqnarray*}
v\left( \overline{x}\right) -v\left( x^{\lambda }\right) &\geq &\left[ \frac{%
\left( n-2\right) }{\sqrt{2}\left\vert x\right\vert ^{n-1}}+O\left( \frac{1}{%
\left\vert x\right\vert ^{n}}\right) \right] \left( y_{1}+x_{1}\right) \\
&=&2\lambda \left[ \frac{\left( n-2\right) }{\sqrt{2}\left\vert x\right\vert
^{n-1}}+O\left( \frac{1}{\left\vert x\right\vert ^{n}}\right) \right] .
\end{eqnarray*}%
It follows that 
\begin{eqnarray*}
v\left( x\right) -v\left( x^{\lambda }\right) &=&v\left( x\right) -v\left( 
\overline{x}\right) +v\left( \overline{x}\right) -v\left( x^{\lambda }\right)
\\
&\geq &0,\text{if }\left\vert x\right\vert \text{ is sufficiently
large,depending on }\lambda _{0}.
\end{eqnarray*}%
This finishes the proof of Lemma \ref{start}.

\begin{lemma}
\label{open}\bigskip If $\lambda _{0}\in \Lambda $, then there exist $%
\varepsilon <0$ s.t. $\left( \lambda _{0}-\varepsilon ,\lambda
_{0}+\varepsilon \right) \subset \Lambda $.
\end{lemma}

The function $v^{\ast }$ defined by $v^{\ast }\left( x\right) =v\left(
x^{\lambda _{0}}\right) $ satisfies%
\begin{equation*}
\begin{array}{ccc}
\Delta v^{\ast }=0 & \text{on} & \mathbb{R}_{+}^{n}, \\ 
-\frac{\partial v}{\partial x_{n}}=\alpha \left( \frac{2}{1+\left\vert
x^{\lambda }\right\vert ^{2}}\right) ^{2}v+\left( \frac{2}{1+\left\vert
x^{\lambda }\right\vert ^{2}}\right) ^{\beta }v^{q} & \text{on} & \mathbb{R}%
^{n-1},%
\end{array}%
\end{equation*}%
Therefore the function $w=v-v^{\ast }$ satisfies, as both $\alpha $ and $%
\beta $ are nonnegative 
\begin{equation*}
\begin{array}{ccc}
\Delta w=0 & \text{on} & \mathbb{R}_{+}^{n}, \\ 
-\frac{\partial w}{\partial x_{n}}\geq qw & \text{on} & \mathbb{R}^{n-1},%
\end{array}%
\end{equation*}%
where 
\begin{equation*}
q=\alpha \left( \frac{2}{1+\left\vert x\right\vert ^{2}}\right) ^{2}+\left( 
\frac{2}{1+\left\vert x\right\vert ^{2}}\right) ^{\beta }\frac{v^{q}-\left(
v^{\ast }\right) ^{q}}{v-v^{\ast }}.
\end{equation*}%
By the assumption we have $w\geq 0$ on $\Sigma _{\lambda _{0}}$. Moreover by
the asymptotic expansion of $v$ it is clear that $w$ cannot be identically
zero on $\Sigma _{\lambda _{0}}$. We claim that $\frac{\partial w}{\partial
x_{1}}<0$ everywhere on the half-plane $\left\{ x:x_{1}=\lambda
_{0},x_{n}\geq 0\right\} $. If $x_{n}>0$, this follows from the Hopf Lemma.
When $x_{n}=0$ one can adapt Hopf's argument to get the the same conclusion,
as observed in Ou \cite{Ou}. We elaborate the idea. First observe that, in
view of the boundary condition for $w$ and the Hopf lemma again, we have $%
w\left( x\right) >0$ on $\left\{ x:x_{1}<\lambda _{0},x_{n}=0\right\} $.
Given $\xi $ with $\xi _{1}=\lambda _{0},\xi _{n}=0$ we consider on the
annulus $A=\left\{ x\in \overline{\mathbb{R}_{+}^{n}}:\rho /2\leq \left\vert
x-\xi -\rho e_{1}\right\vert \leq \rho \right\} $ the following function 
\begin{equation*}
\widetilde{w}\left( x\right) =w\left( x\right) +\delta \left( e^{-\alpha
\rho ^{2}}-e^{-\alpha \left\vert x-\xi -\rho e_{1}\right\vert ^{2}}\right) .
\end{equation*}%
Note that $\widetilde{w}=0$ on $\partial A\cap \left\{ \left\vert x-\xi
-\rho e_{1}\right\vert =\rho \right\} $. For appropriately chosen positive
numbers $\rho ,\alpha ,$ and $\delta $, we can ganrantee%
\begin{eqnarray*}
\Delta \widetilde{w} &\leq &0\text{ on }A; \\
\widetilde{w} &\geq &0\text{ on }\partial A\cap \left\{ \left\vert x-\xi
-\rho e_{1}\right\vert =\rho /2\right\}
\end{eqnarray*}%
while on $\partial A\cap \mathbb{R}^{n-1}$ we have $-\frac{\partial 
\widetilde{w}}{\partial x_{n}}=-\frac{\partial w}{\partial x_{n}}\geq 0$. By
the maximum principle $\widetilde{w}\geq 0$ on $A$. As $\widetilde{w}\left(
\xi \right) =0$ we must have $\frac{\partial \widetilde{w}}{\partial x_{1}}%
\left( \xi \right) \leq 0$. Hence $\frac{\partial w}{\partial x_{1}}\left(
\xi \right) =\frac{\partial \widetilde{w}}{\partial x_{1}}\left( \xi \right)
-\delta \alpha <0$.

We now finish the proof of Lemma \ref{open} by contradiction. Suppose there
is a sequence of positive $\mu _{k}\rightarrow \lambda _{0}$ and $x_{k}\in
\Sigma _{\mu _{k}}$ s.t. $v\left( x_{k}\right) <v\left( x_{k}^{\mu
_{k}}\right) $. By Lemma \ref{start} $\left\{ x_{k}\right\} $ is bounded.
Passing to a subsequence we assume $x_{k}\rightarrow \overline{x}$ as $%
k\rightarrow \infty $. Then $v\left( \overline{x}\right) \leq v\left( 
\overline{x}^{\lambda _{0}}\right) $. This implies that $\overline{x}$ lies
in the half-plane $\left\{ x:x_{1}=\lambda _{0},x_{n}\geq 0\right\} $ and
hence $\frac{\partial u}{\partial x_{1}}\left( \overline{x}\right) <0$. On
the other hand by the mean vaule theorem there exists $\xi _{k}$ in the
segement joining $x_{k}$ and $x_{k}^{\mu _{k}}$ s.t. $\frac{\partial u}{%
\partial x_{1}}\left( \xi _{k}\right) >0$. In the limit we obtain $\frac{%
\partial u}{\partial x_{1}}\left( \overline{x}\right) \geq 0$, a
contradiction.

By Lemma \ref{open} $\Lambda $ is open in $\left( 0,\infty \right) $. It is
also clearly closed. Therefore $\Lambda =\left( 0,\infty \right) $ and we
have $v\left( x_{1},x^{\prime }\right) \geq v\left( -x_{1},x^{\prime
}\right) $ when $x_{1}\leq 0$. As the equation is symmetric w.r.t. $\left(
x_{1},x^{\prime }\right) \rightarrow \left( -x_{1},x^{\prime }\right) $, we
must have $v\left( x_{1},x^{\prime }\right) =v\left( -x_{1},x^{\prime
}\right) $, i.e. $v$ is even in $x_{1}$.

\section{Back to $\overline{\mathbb{B}^{n}}$}

Before we finish the proof of Theorem \ref{main}, we need to explain another
ingredient. Let $(M^{n},g)$ be a compact Riemannian manifold with bounadary $%
\Sigma $ which may be empty. We denote by $T$ the Einstein tensor, i.e. $%
T=Ric-\frac{R}{n}g$, where $R$ is the scalar curvature. If $\Xi $ is a
conformal vector field, then according to \cite{S} the following identity
holds 
\begin{equation*}
\int_{M}\Xi Rdv_{g}=\frac{2n}{n-2}\int_{\Sigma }T(\Xi ,\nu )d\sigma _{g}.
\end{equation*}

We further assume that the boundary is umbilic, i.e. the 2nd fundmental fom
is a proportional to the 1st fundamental form. More precisely for any $%
X,Y\in T\Sigma $%
\begin{equation*}
\Pi \left( X,Y\right) =\frac{H}{n-1}g\left( X,Y\right) ,
\end{equation*}%
where $H$ denotes the mean curvature.

By the Codazzi equation we have for any $X,Y,Z\in T\Sigma $

\begin{eqnarray*}
R(X,Y,Z,\nu ) &=&\nabla _{X}\Pi \left( Y,Z\right) -\nabla _{Y}\Pi \left(
X,Z\right) \\
&=&\frac{XH}{n-1}g\left( Y,Z\right) -\frac{YH}{n-1}g\left( X,Z\right) .
\end{eqnarray*}%
Taking trace yields 
\begin{equation*}
T\left( X,\nu \right) =Ric\left( X,\nu \right) =-\frac{n-2}{n-1}XH.
\end{equation*}%
Therefore we obtain

\begin{proposition}
\label{KW}\bigskip Suppose $(M^{n},g)$ is compact with an umbilic boundary
and $\Xi $ a conformal vector field on $M$ s.t. $\Xi $ is tangential on the
boundary, then%
\begin{equation*}
\int_{M}\Xi Rdv_{g}=-\frac{2n}{n-1}\int_{\Sigma }\Xi Hd\sigma _{g}.
\end{equation*}
\end{proposition}

When $u>0$ satisfies (\ref{E}) the metric $g=u^{4/\left( n-2\right) }dx^{2}$
on $\overline{\mathbb{B}^{n}}$ has scalar curvature $R=0$ and mean curvature%
\begin{equation*}
H=\left[ \left( \frac{n-2}{2}-a\right) u+u^{q}\right] u^{-n/\left(
n-2\right) }=\left( \frac{n-2}{2}-a\right) u^{-2/\left( n-2\right)
}+u^{q-n/\left( n-2\right) }.
\end{equation*}%
on the boundary $\mathbb{S}^{n-1}$. Being umbilic is conformally invariant.
Therefore we can apply Proposition \ref{KW} in this situation. For each $%
i=1,\cdots ,n$ 
\begin{equation*}
\Xi _{i}\left( x\right) =x_{i}x-\frac{1+\left\vert x\right\vert ^{2}}{2}e_{i}
\end{equation*}%
is a conformal vector field on $\overline{\mathbb{B}^{n}}$ and its
restriction on the boundary $\mathbb{S}^{n-1}$ is given by%
\begin{equation*}
\Xi _{i}\left( \xi \right) =\xi _{i}\xi -e_{i}=\nabla \xi _{i}.
\end{equation*}%
Therefore by Propostion \ref{KW} we have for $i=1,\cdots ,n$ 
\begin{equation*}
\int_{\mathbb{S}^{n-1}}\left\langle \nabla \left[ \left( \frac{n-2}{2}%
-a\right) u^{--2/\left( n-2\right) }+u^{q-n/\left( n-2\right) }\right]
,\nabla \xi _{i}\right\rangle u^{2\left( n-1\right) /\left( n-2\right)
}d\sigma =0,
\end{equation*}%
here the gradient, the pairing and the volume element $d\sigma $ are all
with respect to the standard metric on $\mathbb{S}^{n-1}$. Simplifying yields%
\begin{equation}
\int_{\mathbb{S}^{n-1}}\left[ \left( \frac{n-2}{2}-a\right) \frac{2}{n-2}%
u+\left( \frac{n}{n-2}-q\right) u^{q}\right] \left\langle \nabla u,\nabla
\xi _{i}\right\rangle d\sigma =0.  \label{kws}
\end{equation}

By Theorem, we know that $u|_{\mathbb{S}^{n-1}}$ is axial symmetric w.r.t.
the $x_{n}$-axis. Thus we write $u\left( \xi \right) =f\left( \xi
_{n}\right) $, with $f$ a smooth function on $\left[ -1,1\right] $. If $f$
has a critical point $t_{0}\in \left( -1,1\right) $, then every point $\xi
\in \mathbb{S}^{n-1}$ with $\xi _{n}=t_{0}$ is a critical point of $u|_{%
\mathbb{S}^{n-1}}$. By theorem, then $u|_{\mathbb{S}^{n-1}}$ has is axial
symmetric w.r.t. the line passing through $0$ and $\xi $. In other words, $%
u\left( x\right) $ only depends on the distance between $x$ and $\xi $. It
is easy to see that then $f$ must then be constant.

If $f$ has no critical point in $\left( -1,1\right) $, then $f^{\prime }$ is
either everywhere positive or everywhere negative. This implies that $%
\left\langle \nabla u,\nabla \xi _{n}\right\rangle $ is either everywhere
positive or everywhere negative. We then have a contradiction with (\ref{kws}%
) when $i=n$.

Therefore $u|_{\mathbb{S}^{n-1}}$ and hence $u$ itself is constant.

\end{document}